\documentclass[reqno,12pt]{amsart}
\usepackage{geometry}
\usepackage{amssymb}
\usepackage[usenames, dvipsnames]{color}
\usepackage{txfonts}
\usepackage{hyperref}
\usepackage{verbatim}
\usepackage{marginnote}
\usepackage{todonotes}
\usepackage[normalem]{ulem}
\usepackage{cancel}
\usepackage{soul}

\theoremstyle{plain}
\newtheorem{theorem}{Theorem}[section]
\newtheorem{lemma}[theorem]{Lemma}
\newtheorem{corollary}[theorem]{Corollary}
\newtheorem{proposition}[theorem]{Proposition}

\theoremstyle{definition}
\newtheorem{definition}[theorem]{Definition}

\theoremstyle{remark}
\newtheorem{remark}[theorem]{Remark}

\newcommand{\bR}{{\mathbb R}}
\newcommand{\bN}{{\mathbb N}}
\newcommand{\bZ}{{\mathbb Z}}

\def\la{\lambda}
\def\La{\Lambda}
\def\t{\tilde}
\def\q{\quad}

\def\th{\theta}
\def\g{\gamma}

\def\Dl{\Delta}
\def\lt{\left}

\def\rt{\right}

\def\i{\infty}
\def\e{\epsilon}

\def\p{\partial}
\def\f{\frac}
\def\na{\nabla}
\def\al{\alpha}

\def\o{\omega}

\def\s{\sqrt}

\def\nn{\nonumber}
\def\be{\begin{equation}}
\def\ee{\end{equation}}
\def\bes{\begin{equation*}}
\def\ees{\end{equation*}}
\def\bali{\begin{aligned}}
\def\eali{\end{aligned}}

\def\pf{\noindent {\bf Proof. \hspace{2mm}}}
\def\bD{{\mathbb D}}
\def\bW{{\mathbb W}}

\def\dD{\dot{\Delta}}

\numberwithin{equation}{section}

\begin{document}
\title[compressible Oldroyd-B model without retardation]{Globally analytical solutions of the compressible Oldroyd-B model without retardation}

\author[X. Pan]{Xinghong Pan}
\address[X. Pan]{School of Mathematics, Nanjing University of Aeronautics and Astronautics, Nanjing 211106, China}

\email{xinghong\_87@nuaa.edu.cn}

\thanks{X. Pan is supported by National Natural Science Foundation of China (No. 11801268, 12031006.)}



\begin{abstract}
In this paper, we prove the global existence of analytical solutions to the compressible Oldroyd-B model without retardation near a non-vacuum equilibrium in $\bR^n$ $(n=2,3)$.  Zero retardation results in zero dissipation in the velocity equation, which is the main difficulty that prevents us to obtain the long time well-posedness of solutions. Through dedicated analysis, we find that the linearized equations of this model have damping effects, which ensures the global-in-time existence of small data solutions. However, the nonlinear quadratic terms have one more order derivative than the linear part and no good structure is discovered to overcome this derivative loss problem. So we can only build the result in the analytical energy space rather than Sobolev space with finite order derivatives.

\end{abstract}
%\today
\maketitle

%\tableofcontents

\section{Introduction}

\q\ The compressible Oldroyd-B model describes the motion of a type of viscoelastic fluid with memory, which is non-Newtonian and the stress tensor of it is not linearly dependent on the deformation tensor. It is governed by following the conservations of mass and momentum, and the constitutive law.
\be\label{CORB}
\lt\{
\bali
&\rho_t+\operatorname{div} (\rho u)=0,\\
&(\rho u)_t+\na\cdot(\rho u\otimes u)-\mu\o(\Dl u+\na \operatorname{div} u)+\nabla\pi=\operatorname{div}\tau,\\
&\nu(\tau_t+u\cdot\na \tau+Q(\tau,\na u))+\tau=2\mu(1-\o)\bD(u).
\eali
\rt.
\ee
where $\rho\in\bR$, $u\in\bR^n\ (n=2,3)$ and $\tau\in \bR^{n\times n} \ (n=2,3)$ (a symmetric matrix) are the density, the velocity and the non-Newtonian part of the stress tensor, respectively. And $\pi\in\bR$ is the scalar pressure. The constants $\mu>0, \nu>0$ and $\o\in[0,1]$ are viscosity, relaxation time, and the coupling constant. The bilinear term $Q$ has the following form:
$$
Q(\tau, \nabla u)=\tau \bW(u)-\bW(u) \tau-b(\bD(u) \tau+\tau \bD(u)),
$$
where $b \in[-1,1]$ is a parameter, $\bD(u)=\frac{1}{2}\left(\nabla u+ \lt(\nabla u\rt)^T\right)$ is the deformation tensor, and $\bW(u)=\frac{1}{2}\left(\nabla u-\lt(\nabla u\rt)^T\right)$ is the vorticity tensor.

 In this paper, we consider the case $\o=0$ which corresponds to the situation of zero retardation as depicted in Oldroyd \cite{Oldroyd:1958PRSL}. For presentation of the physical background of this compressible model without retardation, we give a brief derivation.  Following \cite{CM:2001SIAM}, an compressible fluid is subject to the following equations
\be\label{incomfluid}
\left\{\begin{aligned}
&\rho_t+\operatorname{div} (\rho u)=0,\\
&(\rho u)_t+\na\cdot(\rho u\otimes u) =\nabla \cdot \sigma, \\
\end{aligned}\right.
\ee
where $\sigma$, an symmetric matrix, is the stress tensor, which can be decomposed as $\sigma=\t{\tau}-\pi \text{\rm Id}$, where $\t{\tau}$ is the tangential part of the stress tensor and $-\pi$Id is the normal part. For a Newtonian fluid, $\t{\tau}$ depends linearly on $\nabla u$ and more precisely
$$
\t{\tau}=2 \mu \bD(u).
$$
However, when we consider the non-Newtonian fluid of Oldroyd-B model, the constitutive law satisfied by $\t{\tau}$ is
\be\label{conlaw}
\t{\tau}+\nu\f{\mathcal{D}\t{\tau}}{\mathcal{D}t}=2\mu\lt(\bD(u)+\t{\nu}\f{\mathcal{D}\bD(u)}{\mathcal{D}t}\rt),
\ee
where for a tensor $\t{\tau}$,
\be\label{1.4}
\f{\mathcal{D}\t{\tau}}{\mathcal{D}t}=\p_t\t{\tau}+u\cdot\na \t{\tau}+Q(\t{\tau},\na u).
\ee
In \eqref{conlaw}, $\nu$ is the relaxation time, $\t{\nu}$ is the retardation time $\left(0 \leq \t{\nu} \leq \nu\right)$, $\mu>0$ is the dynamical viscosity of the fluid. Fluids of this type have both elastic properties and viscous properties. We divide the Oldroyd-B model into the following three types of mathematical system according to the choice of the retardation time and the relaxation time.

Case $0=\t{\nu}=\nu$: This corresponds to purely viscous case (compressible Navier-Stokes equation). This model can be founded in \cite{Lamb:1932CAMB,LandauL:1959} and rigorous mathematical study on the compressible Navier-Stokes equations was initiated since then. By now, there are huge literatures concerning on well-posedness of the compressible Navier-stokes equations. Since in this paper, we consider the non-Newtonian fluid of Oldroyd-B model, we do not pursue any further references details. Readers can refer to citations of\cite{Lamb:1932CAMB,LandauL:1959} to seek more detailed study on the compressible Navier-stokes equations.

Case $0<\t{\nu}\leq \nu$:  If we define
\bes
\tau:=\t{\tau}-2\mu\f{\t{\nu}}{\nu}\bD(u),
\ees
then from \eqref{conlaw}, the second equation of \eqref{incomfluid} and \eqref{1.4}, we see that
\bes
\lt\{
\bali
&(\rho u)_t+\na\cdot(\rho u\otimes u)-\mu\f{\t{\nu}}{\nu}(\Dl u+\na \operatorname{div} u)+\nabla\pi=\operatorname{div}\tau\\
&\nu(\tau_t+u\cdot\na \tau+Q(\tau,\na u))+\tau=2\mu(1-\f{\t{\nu}}{\nu})\bD(u).
\eali
\rt.
\ees
 This corresponds to system \eqref{CORB} with $\o=\f{\t{\nu}}{\nu}\in (0,1]$. This model is an extension of the classical incompressible Oldroyd-B system introduced by Oldroyld in \cite{Oldroyd:1958PRSL}, which has been studied extensively in the literature at present.

Case $0=\t{\nu}<\nu$: This is the purely elastic case (the Maxwell model), which corresponds to \eqref{CORB} with $\o=0$. As far as the author knows, there is little literature to consider this situation. In this paper, we consider this case and the domain is chosen to $\bR^3$.  The exact values of the positive constants $\mu$ and $\nu$ play no essential roles in our following analysis. We set $\mu=\nu=1$ for simplicity. Then system \eqref{CORB} becomes

\be\label{cobr3}
\lt\{
\begin{array}{ll}
\bali
&\rho_t+\operatorname{div}(\rho u)=0,\\
&(\rho u)_t+\na\cdot(\rho u\otimes u)+\nabla\pi=\operatorname{div}\tau,\\
&\tau_t+u\cdot\na \tau+\tau+Q(\tau,\na u)=2\bD(u),
\eali
& (x,t)\in \bR^n\times(0,+\i),\\
(\rho,u,\tau)\big|_{t=0}=(\rho_0, u_0,\tau_0), & x\in\bR^n.
\end{array}
\rt.
\ee
 We will consider that $(\rho,u,\tau)$ is a small analytical perturbation around the non-vacuum equilibrium $(1,0,0)$. Also for simplicity, we consider the pressure $\pi$ satisfying the $\g-$law and assume that $\pi(\rho)=\f{1}{\g}\rho^\g$.

We show the global existence of analytic solutions of system \eqref{cobr3} in Besov spaces by using Fourier analysis method. We have the following theorem.

\begin{theorem}\label{thr3}
Let $\la_0>0$ be a fixed constant. Assume that the initial data $(\rho_0,u_0,\tau_0)$ satisfy $e^{2\la_0\La}(\rho_0-1,u_0,\tau_0)\in \dot{B}^{\f{n}{2}-1}\cap \dot{B}^{\f{n}{2}}$. Then there exists two positive constants $\e_0$ and $C_0$, independent of $\la_0$, such that for any $0<\e\leq\e_0$, if
\bes
\lt\|e^{2\la_0\La}(\rho_0-1,u_0,\tau_0)\rt\|_{\dot{B}^{\f{n}{2}-1}\cap \dot{B}^{\f{n}{2}}}\leq \e\la_0,
\ees
then system \eqref{cobr3} admits a unique solution $(\rho(t),u(t),\tau(t))$ satisfying for any $t\in (0,+\i)$,
\bes
\lt\|e^{\la_0\La}(\rho-1,u,\tau)(t)\rt\|_{\dot{B}^{\f{n}{2}-1}\cap \dot{B}^{\f{n}{2}}}+\lt\|e^{\la_0\La}(\rho-1,u,\tau)\rt\|_{L^1_t(\dot{B}^{\f{n}{2}})}\leq C_0\lt\|e^{2\la_0\La}(\rho_0-1,u_0,\tau_0)\rt\|_{\dot{B}^{\f{n}{2}-1}\cap \dot{B}^{\f{n}{2}}}.
\ees
\end{theorem}

\begin{remark}
Notations of Besov spaces in Theorem \ref{thr3} can be found in Section \ref{sec2}.
\end{remark}

 Before ending of the introduction, we give some literature review related to fluid  model of Oldroyd-B type. For the incompressible Oldroyd-B model, Guillop\'{e}-Saut \cite{GuillopeS:1990NATMA} gave the local well-posedness of regular solution and global existence of small smooth solutions in a smooth open domain with no slip boundary. Lions-Masmoudi \cite{LionsM:2000CAMB} constructed global weak solutions for general initial conditions with the assumptions $b=0$. Chemin-Masmoudi \cite{CM:2001SIAM} gave the existence and uniqueness of locally large and globally small solutions in the critical Besov space. Some remarks on the blow up criteria are shown in Lei-Masmoudi-Zhou \cite{LeiMZ:2010JDE}. Readers can also refer to \cite{ChenM:2008NATMA,ZiFZ:2014ARMA,FangZ:2016SIAM} for more global existence results in Besov spaces. If we consider the case  $\o=0$ but the equation of $\tau$ contains a diffusion  term $-\Dl\tau$, Elgindi-Rousset \cite{ElgindiR:2015CPAM} proved the global existence of small smooth solutions and a similar result with general data if $Q(\tau, \na u)=0$  in 2D.  See also Elgindi-Liu \cite{ElgindiL:2015JDE} for the 3D result. Zhu \cite{Zhu:2018JFA} proved global existence of small solutions without damping effect ($\tau$ is missing in the constitutive law) in Sobolev spaces. See \cite{ChenH:2019JMFM,Zhai:2021JMP} for extensions results in critical Besov spaces, respectively. Recently, Zi \cite{Zi:2021AIHP} consider the vanishing viscosity limits of 3D incompressible Oldroyd-B model in analytical spaces.

For the compressible Oldroyd-B model, there are relatively fewer literatures.  The incompressible limit problems in torus and bounded domain were investigated in Lei \cite{Lei:2006CAMB}, and Guillop\'{e}-Salloum-Talhouk \cite{GuillopeST:2010DCDSB}, respectively for well-prepared data. The case of ill-prepared initial data was studied by Fang and Zi \cite{FangZ:2014JDE}.  Global well-posedness and decay rates results in $H^2$-framework for the three dimensional case was given in Zhou-Zh-Zi \cite{ZhouZZ:2018JDE}.  See Zhu \cite{Zhu:2022CVPDE} and Pan-Xu-Zhu \cite{PanXZ:2022JDE} for global existence results in Sobolev spaces and critical Besov spaces for the model without damping mechanism. Readers can refer to \cite{LeiZ:2005SIAM,QianZ:2010ARMA,HuW:2011JDE,PanX:2019DCDS} and references therein for more compressible Oldroyd-type model results.

Throughout the paper, $C_{a,b,c,...}$ denotes a positive constant depending on $a,\,b,\, c,\,...$ which may be different from line to line. We also apply $A\lesssim_{a,b,c,\cdots} B$ to denote $A\leq C_{a,b,c,...}B$. $A\thickapprox_{a,b,c,\cdots} B$ means $A\lesssim_{a,b,c,\cdots} B$ and $B\lesssim_{a,b,c,\cdots}A$. For a norm $\|\cdot\|$, we use $\|(f,g,\cdots)\|$ to denote $\|f\|+\|g\|+\cdots$. For a function $f(t,x)$, $\|f(t)\|_{L^p(\bR^n)}$ denote the usual spacial $L^p$ norm for $1\leq p\leq +\i$. Besides, if $p=2$, we will simply denote $\|f(t)\|_{L^2(\bR^n)}$ by $\|f(t)\|$. We use $[A, B]=AB-BA$ to denote the commutator of $A$ and $B$.

Our paper is arranged as follows. In Section \ref{sec2}, we give an a priori estimate solution of system \eqref{cobr3}. Then The proof of Theorem \ref{thr3} is obtained in Section \ref{sec3}.

\section{The a priori estimate}\label{sec2}

In this section, we give an a priori estimate for system \eqref{cobr3}. First we introduce some notations and the functional spaces which we use.

\subsubsection*{Notations}
\begin{itemize}
\item For a function $f\in\mathcal{S}'$ (the dual space of the Schwartz space), denote by $\hat{f}$ or $\mathcal{F}(f)$ the Fourier transform of $f$, and $\check{f}$ or $\mathcal{F}^{-1}(f)$ the inverse Fourier transform of $f$;
\item Denote $\s{-\Dl}$ by $\La$ and for any $s\in\bR$, $\La^s=(-\Dl)^{\f{s}{2}}$;
\item For functions $f,g\in L^2$, denote the $L^2$ inner product by $(f|g)$, namely,
\bes
(f|g)=\int_{\bR^3} f\bar{g}dx;
\ees
\item Denote by\ $\mathcal{Z}'(\bR^n)$ the dual space of
\bes
\mathcal{Z}(\bR^n):=\{f\in \mathcal{S}(\bR^n):\p^\al \hat{f}(0)=0,\forall \al\in(\bN\cup 0)^n\}.
\ees
\end{itemize}

\subsubsection*{Littlewood-Paley decomposition}

Next, we need a Littlewood-Paley decomposition. There exist two radial smooth functions $\varphi(x),\,\chi(x)$ supported in the annulus $\mathcal{C}=\{\xi\in\bR^n:3/4\leq |\xi|\leq 8/3\}$ and the ball $B=\{\xi\in\bR^n:|\xi|\leq 4/3\}$, respectively such that
\be
\sum\limits_{j\in \bZ}\varphi(2^{-j}\xi)=1\q \forall \xi\in\bR^n\setminus\{0\}.\nn
\ee
The homogeneous dyadic blocks $\dD_j$ and the homogeneous low-frequency cut-off operators $\dot{S}_j$ are defined for all $j\in \bZ$ by
\be
\dD_j u=\varphi(2^{-j}D)f,\q \dot{S}_jf=\sum\limits_{k\leq j-1}\dD_k f=\chi(2^{-j}D)f.\nn
\ee

Let us now turn to the definition of the main functional spaces and norms that
will come into play in our paper.
\begin{definition}
Let $s$ be a real number and $(p,r)$ be in $[1,\i]^2$. The homogeneous Besov space $\dot{B}^s_{p,r}$ consists of those distributions $u\in \mathcal{Z}'(\bR^n)$ such that
\be
\|u\|_{\dot{B}^s_{p,r}}\triangleq \Big(\sum\limits_{j\in\bZ}2^{jsr}\|\dD_j u\|^r_{L^p}\Big)^{\f{1}{r}}<\i.\nn
\ee
\end{definition}
Also, we introduce the hybrid Besov space since our analysis will be performed at different frequencies.
\begin{definition}
Let $s,\sigma\in \bR$. The hybrid Besov space $\dot{B}^{s,\sigma}$ is defined by
\be
\dot{B}^{s,\sigma}\triangleq \{f\in\mathcal{Z}'(\bR^n):\|f\|_{\dot{B}^{s,\sigma}}<\i\}, \nn
\ee
with
\begin{align}
\|f\|_{\dot{B}^{s,\sigma}}\triangleq& \sum\limits_{k\leq k_0}2^{ks}\|\dD_k f \|_{L^2}+\sum\limits_{k>k_0}2^{k\sigma}\|\dD_k f \|_{L^2},\nn\\
                          =&\|f\|^\ell_{\dot{B}^s_{2,1}}+\|f\|^h_{\dot{B}^\sigma_{2,1}},\nn
\end{align}
\end{definition}
where $k_0$ is a fixed suitably large constant to be defined.
\begin{remark} We note that
\begin{itemize}
\item If $\sigma=s$, $\dot{B}^{s,s}$ is the usual Besov space $\dot{B}^s_{2,1}$;
\item If $\sigma<s$, $\dot{B}^{s,\sigma}=\dot{B}^s_{2,1}\cap\dot{B}^\sigma_{2,1}$.
\end{itemize}
\end{remark}
In the case where $u$ depends on the time variable, we consider the space-time mixed spaces as follows
\be
\|u\|_{L^q_T\dot{B}^{s,\sigma}}:=\big\|\|u(t,\cdot)\|_{\dot{B}^{s,\sigma}}\big\|_{L^q(0,T)}.\nn
\ee
In addition, we introduce another space-time mixed spaces, which is usually referred to Chemin-Lerner's spaces. The definition is given by
\be
\|u\|_{\t{L}^q_T\dot{B}^{s,\sigma}}\triangleq\sum\limits_{k\leq k_0}2^{ks}\|\dD_k u\|_{L^q(0,T)L^2}+\sum\limits_{k> k_0}2^{k\sigma}\|\dD_k u\|_{L^q(0,T)L^2}.\nn
\ee
The index $T$ will be omitted if $T=+\i$. It is easy to check that $\t{L}^1_T\dot{B}^{s,\sigma}=L^1_T\dot{B}^{s,\sigma}$ and $\t{L}^q_T\dot{B}^{s,\sigma}\subseteq L^q_T\dot{B}^{s,\sigma}$ for $q>1$.

In this paper, we also need the following time-weighted hybrid Besov norm.
\begin{definition}
Let $\th(t)\in L^1_{\operatorname{loc}}(\bR_+)$ be a positive function. Define
\bes
\|f\|_{L^1_{T,\th(t)}(\dot{B}^{s,\sigma})}=\int^T_0\th(t)\|f(t)\|_{\dot{B}^{s,\sigma}}dt.
\ees
\end{definition}
Let $\th(t)\in C^1[0,+\i)$ with $\th(0)=0$ be a non-decreasing function. Denote
\bes
\Phi(t,\xi)=(2\la_0-\la\th(t))|\xi|,
\ees
where $\la$ is suitably large constant and will be determined later. For a function $f$, define
\be\label{fourier}
f_{\Phi}(t,x)= \mathcal{F}^{-1}_{\xi\rightarrow x}\lt(e^{\Phi(t,\xi)}\hat{f}(t,\xi)\rt)=e^{\Phi(t,\La)}{f}(t,x).
\ee
In particular,
\bes
f_{\Phi}(0,x)=e^{2\la_0\La}{f}(0,x),\,\ \text{and }\ \dD_k f_{\Phi}(t,x)= \mathcal{F}^{-1}_{\xi\rightarrow x}\lt(\varphi(2^{-k}\xi)e^{\Phi(t,\xi)}\hat{f}(t,\xi)\rt).
\ees
Later for convenience and simplification of notations, we use $f_{\Phi,k}$ to denote $\dD_k f_{\Phi}$.

Now our a priori estimate is stated in the following Proposition.
\begin{proposition}
Let $\la_0>0$ be a fixed constant. Assume that $(\rho,u,\tau)$ is a solution of of system \ref{cobr3}, with the initial data $(\rho_0,u_0,\tau_0)$ satisfying $e^{2\la_0\La}(\rho_0-1,u_0,\tau_0)\in \dot{B}^{n/2-1,n/2}$. Then there exists a uniform constant $C$, such that for any $t>0$, we have
\be\label{apriori}
\bali
&\|(a_{\Phi},u_{\Phi},\tau_{\Phi})(t)\|_{\dot{B}^{n/2-1,n/2}}+\la \|(a_{\Phi},u_{\Phi},\tau_{\Phi})\|_{L^1_{t,\dot{\th}}(\dot{B}^{n/2,n/2+1})}+\|(a_{\Phi},u_{\Phi},\tau_{\Phi})\|_{L^1_{t}(\dot{B}^{n/2}_{2,1})}\\
\leq& C\lt(\|e^{2\la_0\La}(a_0,u_0,\tau_0)\|_{\dot{B}^{n/2-1,n/2}}+\int^t_0\|(a_{\Phi},u_{\Phi},\tau_{\Phi})(s)\|_{\dot{B}^{n/2}_{2,1}}\|( a_{\Phi}, u_{\Phi}, \tau_{\Phi})(s)\|_{{B}^{n/2,n/2+1}_{2,1}}ds\rt).
\eali
\ee
Here $a:=\rho-1$, $a_0:=\rho_0-1$, and the constant $C$ is independent of $\la_0$.
\end{proposition}

\subsection{Linearized problem and its estimates}

In this section, we linearize system \eqref{cobr3} and give its linear a priori estimates.

Set $\rho=1+a$, then we rewrite \eqref{cobr3} into
\be\label{cobr3first}
\lt\{
\begin{array}{ll}
\bali
&a_t+\operatorname{div}u =F,\\
&u_t+\na a-\operatorname{div}\tau=G,\\
&\tau_t+\tau-2\bD(u)=H,
\eali
& (x,t)\in \bR^3\times(0,+\i),\\
(a,u,\tau)\big|_{t=0}=(a_0, u_0,\tau_0), & x\in\bR^3,
\end{array}
\rt.
\ee
where
\begin{align}
&F:=-\na\cdot(au),\nn\\
&G:=-u\cdot \na u+[1-(1+a)^{\g-2}]\na a-\f{a}{1+a}\operatorname{div}\tau,\label{fgh}\\
&H:=-u\cdot\na \tau-Q(\tau,\na u).\nn
\end{align}

Then, we have the following a priori estimate for the linearized system \eqref{cobr3first}. .
\begin{proposition}
Assume that $(a,u,\tau)$ is a solution of of the linearized system \ref{cobr3first},with the initial data $(a_0,u_0,\tau_0)$ satisfying $e^{2\la_0\La}(a_0,u_0,\tau_0)\in \dot{B}^{n/2-1,n/2}$. Then there exists a uniform constant $C$, such that for any $t>0$, we have
\be\label{apriori}
\bali
&\|(a_{\Phi},u_{\Phi},\tau_{\Phi})(t)\|_{\dot{B}^{n/2-1,n/2}}+\la \|(a_{\Phi},u_{\Phi},\tau_{\Phi})\|_{L^1_{t,\dot{\th}}(\dot{B}^{n/2,n/2+1})}+\|(a_{\Phi},u_{\Phi},\tau_{\Phi})\|_{L^1_{t}(\dot{B}^{n/2}_{2,1})}\\
\leq& C\lt(\|e^{2\la_0\La}(a_0,u_0,\tau_0)\|_{\dot{B}^{n/2-1,n/2}}+\|(F_{\Phi},G_{\Phi},H_{\Phi})\|_{L^1_{t}(\dot{B}^{n/2-1,n/2})}\rt).
\eali
\ee
\end{proposition}
\pf From \eqref{cobr3first}, we see that $(a_{\Phi}, u_{\Phi}, \tau_{\Phi})$ satisfy the following equations.

\be\label{cobr3second}
\lt\{
\begin{array}{ll}
\bali
&\p_t a_{\Phi,k}+\la \dot{\th}(t) \La a_{\Phi,k}+\operatorname{div}u_{\Phi,k} =F_{\Phi,k},\\
& \p_t u_{\Phi,k}+\la \dot{\th}(t) \La u_{\Phi,k}+\na a_{\Phi,k}-\operatorname{div}\tau_{\Phi,k}=G_{\Phi,k},\\
&\p_t \tau_{\Phi,k}+\la \dot{\th}(t) \La \tau_{\Phi,k}+\tau_{\Phi,k}-2\bD u_{\Phi,k}=H_{\Phi,k}.
\eali
\end{array}
\rt.
\ee
Performing $L^2$ inner product of \eqref{cobr3second}$_{1,2,3}$ with $(2a_{\Phi,k},2u_{\Phi,k},\tau_{\Phi,k})$ respectively, we can obtain that
\begin{align}
&\f{1}{2}\f{d}{dt}\lt(2\|a_{\Phi,k}\|^2_{L^2}+2\|u_{\Phi,k}\|^2_{L^2}+\|\tau_{\Phi,k}\|^2_{L^2}\rt)\nn\\
&+\la\dot{\th}(t)\lt(2\|\La^{1/2}a_{\Phi,k}\|^2_{L^2}+2\|\La^{1/2}u_{\Phi,k}\|^2_{L^2}+\|\La^{1/2}\tau_{\Phi,k}\|^2_{L^2}\rt)+\|\tau_{\Phi,k}\|^2_{L^2}\nn\\
&+2(a_{\Phi,k}|\text{div} u_{\Phi,k})+2(\na a_{\Phi,k}| u_{\Phi,k})-2(\text{div}\tau_{\Phi,k}| u_{\Phi,k})-2(\tau_{\Phi,k}|\bD u_{\Phi,k})\nn\\
=&2(a_{\Phi,k}|F_{\Phi,k})+2(G_{\Phi,k}| u_{\Phi,k})+(\tau_{\Phi,k}|H_{\Phi,k}).\nn
\end{align}
By using the symmetry of $\tau$ and integration by parts, we see that terms on the third line of the above equality are cancelled. Then we can obtain that
\begin{align}
&\f{d}{dt}\lt(2\|a_{\Phi,k}\|^2_{L^2}+2\|u_{\Phi,k}\|^2_{L^2}+\|\tau_{\Phi,k}\|^2_{L^2}\rt)\nn\\
&+2\la\dot{\th}(t)\lt(2\|\La^{1/2}a_{\Phi,k}\|^2_{L^2}+2\|\La^{1/2}u_{\Phi,k}\|^2_{L^2}+\|\La^{1/2}\tau_{\Phi,k}\|^2_{L^2}\rt)+2\|\tau_{\Phi,k}\|^2_{L^2}\label{high0}\\
=&4(a_{\Phi,k}|F_{\Phi,k})+4(G_{\Phi,k}| u_{\Phi,k})+2(\tau_{\Phi,k}|H_{\Phi,k}).\nn
\end{align}
\subsubsection*{Estimate in high frequences}
We apply $\La^{-1}\na$ to \eqref{cobr3second}$_1$, $\La^{-1}$ to \eqref{cobr3second}$_2$ and $\La^{-1}\text{div}$ to \eqref{cobr3second}$_3$ to obtain that
\be\label{cobr3third}
\lt\{
\begin{array}{ll}
\bali
&\p_t \La^{-1}\na a_{\Phi,k}+\la \dot{\th}(t) \na a_{\Phi,k}+\La^{-1}\na \operatorname{div}u_{\Phi,k} =\La^{-1}\na F_{\Phi,k},\\
& \p_t \La^{-1} u_{\Phi,k}+\la \dot{\th}(t) u_{\Phi,k}+\La^{-1} \na a_{\Phi,k}-\La^{-1}\operatorname{div}\tau_{\Phi,k}=\La^{-1}G_{\Phi,k},\\
&\p_t \La^{-1}\text{div}\tau_{\Phi,k}+\la \dot{\th}(t) \text{div} \tau_{\Phi,k}+\La^{-1}\text{div} \tau_{\Phi,k}+(\La-\na  \La^{-1}\text{div})u_{\Phi,k}=\La^{-1}\text{div} H_{\Phi,k}.
\eali
\end{array}
\rt.
\ee
Multiplying \eqref{cobr3third}$_1$ by $\La^{-1} u_{\Phi,k}$ and \eqref{cobr3third}$_2$ by $\La^{-1}\na a_{\Phi,k}$, then integrating over $\bR^n$ to obtain that
\begin{align}
&\p_t(\La^{-1}\na a_{\Phi,k}|\La^{-1} u_{\Phi,k})+2\la \dot{\th}(t) (\na a_{\Phi,k}|\La^{-1} u_{\Phi,k})-\|\La^{-1}\operatorname{div}u_{\Phi,k}\|^2_{L^2}\nn\\
&+\| a_{\Phi,k}\|^2_{L^2}-(\La^{-1}\na a_{\Phi,k}|\La^{-1}\operatorname{div}\tau_{\Phi,k}) \label{high1}\\
=&(\La^{-1}\na F_{\Phi,k}|\La^{-1} u_{\Phi,k})+(\La^{-1}\na a_{\Phi,k}|\La^{-1}G_{\Phi,k}),\nn
\end{align}
where we have used integration by parts to obtain that
\bes
(\La^{-1}\na \operatorname{div}u_{\Phi,k}|\La^{-1} u_{\Phi,k})=-\|\La^{-1}\operatorname{div}u_{\Phi,k}\|^2_{L^2}.
\ees
Multiplying \eqref{cobr3third}$_2$ by $\La^{-1}\text{div}\tau_{\Phi,k}$ and \eqref{cobr3third}$_3$ by $\La^{-1} u_{\Phi,k}$, then integrating over $\bR^n$ to obtain that
\begin{align}
&\p_t(\La^{-1} u_{\Phi,k}|\La^{-1}\text{div}\tau_{\Phi,k})+2\la \dot{\th}(t) (\La^{-1} u_{\Phi,k}|\text{div}\tau_{\Phi,k})+(\La^{-1}\na a_{\Phi,k}|\La^{-1}\operatorname{div}\tau_{\Phi,k})-\| \La^{-1}\text{div}\tau_{\Phi,k}\|^2_{L^2}\nn\\
&+(\La^{-1} u_{\Phi,k}|\La^{-1} \operatorname{div}\tau_{\Phi,k})+\| u_{\Phi,k}\|^2_{L^2}+\| \La^{-1}\text{div}u_{\Phi,k}\|^2_{L^2}\label{high2}\\
=&(\La^{-1}G_{\Phi,k}|\La^{-1}\text{div}\tau_{\Phi,k})+(\La^{-1}u_{\Phi,k}|\La^{-1}\text{div}H_{\Phi,k}). \nn
\end{align}
Adding \eqref{high1} and \eqref{high2} together indicates that
\begin{align}
&\p_t\lt\{(\La^{-1}\na a_{\Phi,k}|\La^{-1} u_{\Phi,k})+(\La^{-1} u_{\Phi,k}|\La^{-1}\text{div}\tau_{\Phi,k})\rt\}\nn\\
&+2\la \dot{\th}(t) \lt\{ (\na a_{\Phi,k}|\La^{-1} u_{\Phi,k})+(\La^{-1} u_{\Phi,k}|\text{div}\tau_{\Phi,k})\rt\}+\| a_{\Phi,k}\|^2_{L^2}+\| u_{\Phi,k}\|^2_{L^2}\nn\\
&+(\La^{-1} u_{\Phi,k}|\La^{-1} \operatorname{div}\tau_{\Phi,k})-\| \La^{-1}\text{div}\tau_{\Phi,k}\|^2_{L^2}\label{high3}\\
=&(\La^{-1}\na F_{\Phi,k}|\La^{-1} u_{\Phi,k})+(\La^{-1}\na a_{\Phi,k}|\La^{-1}G_{\Phi,k})+(\La^{-1}G_{\Phi,k}|\La^{-1}\text{div}\tau_{\Phi,k})+(\La^{-1}u_{\Phi,k}|\La^{-1}\text{div}H_{\Phi,k}). \nn
\end{align}
From the Bernstein's inequality, it is easy to see that
\be\label{high4}
\f{3}{4}2^k \| \dD_kf\|_{L^2}\leq \|\La \dD_kf\|_{L^2}\leq \f{8}{3}2^k\| \dD_kf\|_{L^2}.
\ee
For $k\geq k_0$, define
\bes
\widetilde{\mathcal{E}}^2_k=2\|a_{\Phi,k}\|^2_{L^2}+2\|u_{\Phi,k}\|^2_{L^2}+\|\tau_{\Phi,k}\|^2_{L^2}+(\La^{-1}\na a_{\Phi,k}|\La^{-1} u_{\Phi,k})+(\La^{-1} u_{\Phi,k}|\La^{-1}\text{div}\tau_{\Phi,k}),
\ees
and
\bes
{\mathcal{E}}^2_k=\|a_{\Phi,k}\|^2_{L^2}+\|u_{\Phi,k}\|^2_{L^2}+\|\tau_{\Phi,k}\|^2_{L^2}.
\ees
We choose $k_0=3$, then by using \eqref{high4}, we can see that
\bes
\f{1}{2}{\mathcal{E}}^2_k\leq \widetilde{\mathcal{E}}^2_k \leq 3{\mathcal{E}}^2_k.
\ees
Adding \eqref{high0} and \eqref{high3} together, and using H\"{o}lder inequality, Cauchy inequality and \eqref{high4}, we can obtain that
\begin{align}
&\f{d}{dt}\widetilde{\mathcal{E}}^2_k+\f{1}{4}\la \dot{\th}(t)2^k\widetilde{\mathcal{E}}^2_k +\f{1}{8}\widetilde{\mathcal{E}}^2_k\leq C\|(F_{\Phi,k},G_{\Phi,k}, H_{\Phi,k})\|_{L^2}\widetilde{\mathcal{E}}_k. \nn
\end{align}
where, $C$ is a uniform constant, independent of $\la_0$.  Then from the above inequality, we have
\begin{align}
&\f{d}{dt}\widetilde{\mathcal{E}}_k+\f{1}{8}\la \dot{\th}(t)2^k\widetilde{\mathcal{E}}_k +\f{1}{16}\widetilde{\mathcal{E}}_k\leq C\|(F_{\Phi,k},G_{\Phi,k}, H_{\Phi,k})\|_{L^2}. \label{high6}
\end{align}
Multiplying \eqref{high6} by $2^{\f{n}{2}k}$, integrating the resulted equation from $0$ to $t$ with the time variable, and then summing over $k_0<k\in\bN $, we can achieve that
\begin{align}
&\|(a_{\Phi},u_{\Phi},\tau_{\Phi})(t)\|^h_{\dot{B}^{n/2}_{2,1}}+\la\int^t_0\dot{\th}(\tau)\|(a_{\Phi},u_{\Phi},\tau_{\Phi})(s)\|^h_{\dot{B}^{n/2+1}_{2,1}}ds
+\int^t_0\|(a_{\Phi},u_{\Phi},\tau_{\Phi})(s)\|^h_{\dot{B}^{n/2}_{2,1}}ds\nn\\
\leq&  C\int^t_0\|(F_{\Phi},G_{\Phi},H_{\Phi})(s)\|^h_{\dot{B}^{n/2}_{2,1}}ds.\label{high7}
\end{align}

\subsubsection*{Estimate in low frequences}
We apply $\na$ to \eqref{cobr3second}$_1$ and $\text{div}$ to \eqref{cobr3second}$_3$ to obtain that
\be\label{low0}
\lt\{
\begin{array}{ll}
\bali
&\p_t \na a_{\Phi,k}+\la \dot{\th}(t) \na \La a_{\Phi,k}+\na \operatorname{div}u_{\Phi,k} =\na F_{\Phi,k},\\
& \p_t u_{\Phi,k}+\la \dot{\th}(t) \La u_{\Phi,k}+ \na a_{\Phi,k}-\operatorname{div}\tau_{\Phi,k}=G_{\Phi,k},\\
&\p_t \text{div}\tau_{\Phi,k}+\la \dot{\th}(t) \text{div} \La\tau_{\Phi,k}+\text{div} \tau_{\Phi,k}-(\Dl+\na \text{div})u_{\Phi,k}=\text{div} H_{\Phi,k}.
\eali
\end{array}
\rt.
\ee
Multiplying \eqref{low0}$_1$ by $ u_{\Phi,k}$ and \eqref{low0}$_2$ by $\na a_{\Phi,k}$, then integrating over $\bR^n$ to obtain that
\begin{align}
&\p_t(\na a_{\Phi,k}| u_{\Phi,k})+2\la \dot{\th}(t) (\na \La a_{\Phi,k}| u_{\Phi,k})-\|\operatorname{div}u_{\Phi,k}\|^2_{L^2}\nn\\
&+\| \na a_{\Phi,k}\|^2_{L^2}-(\na a_{\Phi,k}|\operatorname{div}\tau_{\Phi,k})=(\na F_{\Phi,k}| u_{\Phi,k})+(\na a_{\Phi,k}|G_{\Phi,k}),\label{low1}
\end{align}
where we have used integration by parts to obtain that
\bes
(\na \operatorname{div}u_{\Phi,k}| u_{\Phi,k})=-\|\operatorname{div}u_{\Phi,k}\|^2_{L^2}.
\ees
Multiplying \eqref{low0}$_2$ by $\text{div}\tau_{\Phi,k}$ and \eqref{low0}$_3$ by $ u_{\Phi,k}$, then integrating over $\bR^n$ to obtain that
\begin{align}
&\p_t( u_{\Phi,k}|\text{div}\tau_{\Phi,k})+2\la \dot{\th}(t) (\La u_{\Phi,k}|\text{div}\tau_{\Phi,k})+(\na a_{\Phi,k}|\operatorname{div}\tau_{\Phi,k})-\| \text{div}\tau_{\Phi,k}\|^2_{L^2}\nn\\
&+( u_{\Phi,k}| \operatorname{div}\tau_{\Phi,k})+\| \La u_{\Phi,k}\|^2_{L^2}+\|\text{div}u_{\Phi,k}\|^2_{L^2}=(G_{\Phi,k}|\text{div}\tau_{\Phi,k})+(u_{\Phi,k}|\text{div}H_{\Phi,k}). \label{low2}
\end{align}
Adding \eqref{low1} and \eqref{low2} together indicates that
\begin{align}
&\p_t\lt\{(\na a_{\Phi,k}| u_{\Phi,k})+( u_{\Phi,k}|\text{div}\tau_{\Phi,k})\rt\}\nn\\
&+2\la \dot{\th}(t) \lt\{ (\na \La a_{\Phi,k}| u_{\Phi,k})+(\La u_{\Phi,k}|\text{div}\tau_{\Phi,k})\rt\}+\| \La a_{\Phi,k}\|^2_{L^2}+\| \La u_{\Phi,k}\|^2_{L^2}\nn\\
&+( u_{\Phi,k}| \operatorname{div}\tau_{\Phi,k})-\|\text{div}\tau_{\Phi,k}\|^2_{L^2}\label{low3}\\
=&(\na F_{\Phi,k}| u_{\Phi,k})+(\na a_{\Phi,k}|G_{\Phi,k})+(G_{\Phi,k}|\text{div}\tau_{\Phi,k})+(u_{\Phi,k}|\text{div}H_{\Phi,k}). \nn
\end{align}

For $k\leq k_0$, define
\bes
\widetilde{\mathcal{E}}^2_k=2\|a_{\Phi,k}\|^2_{L^2}+2\|u_{\Phi,k}\|^2_{L^2}+\|\tau_{\Phi,k}\|^2_{L^2}+\f{3}{8}2^{-k}(\na a_{\Phi,k}| u_{\Phi,k})+\f{3}{8}2^{-k}(u_{\Phi,k}|\text{div}\tau_{\Phi,k}),
\ees
and
\bes
{\mathcal{E}}^2_k=\|a_{\Phi,k}\|^2_{L^2}+\|u_{\Phi,k}\|^2_{L^2}+\|\tau_{\Phi,k}\|^2_{L^2}.
\ees
We choose $k_0=3$, then by using \eqref{high4}, we can see that
\bes
\f{1}{2}{\mathcal{E}}^2_k\leq \widetilde{\mathcal{E}}^2_k \leq 3{\mathcal{E}}^2_k.
\ees
Multiplying \eqref{low3} by $\f{3}{8}2^{-k}$,  adding the resulted equation to \eqref{high0}, and using H\"{o}lder inequality, Cauchy inequality and \eqref{high4}, we can obtain that
\begin{align}
&\f{d}{dt}\widetilde{\mathcal{E}}^2_k+\f{1}{4}\la \dot{\th}(t)2^k\widetilde{\mathcal{E}}^2_k +\f{1}{8}2^{k}\widetilde{\mathcal{E}}^2_k\leq C\|(F_{\Phi,k},G_{\Phi,k}, H_{\Phi,k})\|_{L^2}\widetilde{\mathcal{E}}_k. \nn
\end{align}
where, $C$ is a uniform constant, independent of $\la_0$.  Then from the above inequality, we have
\begin{align}
&\f{d}{dt}\widetilde{\mathcal{E}}_k+\f{1}{8}\la \dot{\th}(t)2^k\widetilde{\mathcal{E}}_k +\f{1}{16}2^{k}\widetilde{\mathcal{E}}_k\leq C\|(F_{\Phi,k},G_{\Phi,k}, H_{\Phi,k})\|_{L^2}. \label{low6}
\end{align}
Multiplying \eqref{low6} by $2^{(n/2-1)k}$, integrating the resulted equation from $0$ to $t$ with the time variable, and then summing over $\bN\ni k\leq k_0$, we can achieve that
\begin{align}
&\|(a_{\Phi},u_{\Phi},\tau_{\Phi})(t)\|^\ell_{\dot{B}^{n/2-1}_{2,1}}+\la\int^t_0\dot{\th}(s)\|(a_{\Phi},u_{\Phi},\tau_{\Phi})(s)\|^\ell_{\dot{B}^{n/2}_{2,1}}ds
+\int^t_0\|(a_{\Phi},u_{\Phi},\tau_{\Phi})(s)\|^\ell_{\dot{B}^{n/2}_{2,1}}ds\nn\\
\leq&  C\int^t_0\|(F_{\Phi},G_{\Phi},H_{\Phi})(s)\|^\ell_{\dot{B}^{n/2-1}_{2,1}}ds.\label{low7}
\end{align}
\subsubsection*{Proof of the a priori estimate in \eqref{apriori}}

By adding \eqref{high7} and \eqref{low7}, we can achieve the a priori estimate in \eqref{apriori}. \qed

\subsection{Estimates of nonlinear terms}

Using the Bony decomposition, we have the following Lemma.
\begin{lemma}\label{lemproductf}
Let $f_{\Phi}$ be defined in \eqref{fourier}. For $s\in (-n/2,n/2]$, and $(f_{\Phi},g_{\Phi})\in \dot{B}^{n/2}_{2,1}\times\dot{B}^s_{2,1}$, there exists a positive constant $C$, depending on $s$, such that the following product estimate holds.
\be\label{fproduct1}
\|(fg)_{\Phi}\|_{\dot{B}^s_{2,1}}\leq C\|f_{\Phi}\|_{\dot{B}^{n/2}_{2,1}}\|g_{\Phi}\|_{\dot{B}^s_{2,1}}.
\ee
\end{lemma}
Proof of Lemma \ref{lemproductf} is postponed to Appendix.  By Letting $s=n/2$ in \eqref{fproduct1}, we have the following Corollary.
\begin{corollary}
Let $f_{\Phi}$ be defined in \eqref{fourier} and $f_{\Phi}\in \dot{B}^{n/2}_{2,1}$, there exists a positive constant $C$, such that for any $k\in\bN/\{0\}$, the following estimate holds.
\be\label{fproduct6}
\|(f^k)_{\Phi}\|_{\dot{B}^{n/2}_{2,1}}\leq C^k\|f_{\Phi}\|^k_{\dot{B}^{n/2}_{2,1}}.
\ee
\end{corollary}
Achievement of \eqref{fproduct6} is $k-$time use of \eqref{fproduct1} with $s=n/2$. \qed

Now, we use product estimates in \eqref{fproduct1} and \eqref{fproduct6} to give estimates of nonlinear terms $F_{\Phi},\ G_{\Phi}$ and $H_{\Phi}$. We have the following estimate.
\begin{lemma} Let $F,\, G$ and $H$ be defined in \eqref{fgh}. There exists a constant $\e_0$ and $C$ such that if
\bes
\|a(t)\|_{\dot{B}^{n/2}_{2,1}}\leq \e_0<1,
\ees
then for $s\in (-n/2,n/2]$, we have
\be\label{fproduct8}
\|(F_{\Phi},G_{\Phi},H_{\Phi})\|_{\dot{B}^{s}_{2,1}}\leq C\|(a_{\Phi},u_{\Phi},\tau_{\Phi})\|_{\dot{B}^{n/2}_{2,1}}\|(a_{\Phi},u_{\Phi},\tau_{\Phi})\|_{\dot{B}^{s+1}_{2,1}}.
\ee
\end{lemma}
\pf  We first deal with $\dot{B}^{n/2}_{2,1}$ norm of
\bes
1-(1+a)^{\g-2}, \q\text{and}\q \f{a}{1+a}.
\ees

By using Taylor expansion, we see that
\bes
1-(1+a)^{\g-2}=-\sum^\i_{k=1}C_{k,\nu} a^k, \q \text{with}\q C_{k,\nu}=\f{(\g-2)(\g-3)\cdots(\g-2-k+1)}{k!}.
\ees
It is easy to see that there exists a constant $C$ such that $C_{\g,k}\leq C^k$. Using \eqref{fproduct6}, we can obtain that
\begin{align}
\lt\| (1-(1+a)^{\g-2})_{\Phi}\rt\|_{\dot{B}^{n/2}_{2,1}} \leq& \sum^\i_{k=1} C_{k,\nu} \lt\|(a^k)_{\Phi}\rt\|_{\dot{B}^{n/2}_{2,1}}\leq \sum^\i_{k=1} C^k \lt\|a_{\Phi}\rt\|^k_{\dot{B}^{n/2}_{2,1}}\nn\\
                \leq &\f{C \lt\|a_{\Phi}\rt\|_{\dot{B}^{n/2}_{2,1}}}{1-C \lt\|a_{\Phi}\rt\|_{\dot{B}^{n/2}_{2,1}}}\leq C \lt\|a_{\Phi}\rt\|_{\dot{B}^{n/2}_{2,1}},\label{fproduct7}
\end{align}
provided that $C\e_0\leq 1/2$.  The same is true for $a/(1+a)$. From the representation formula of $F, G$ and $H$, by using \eqref{fproduct1} and the above estimate \eqref{fproduct7}, we have that
\begin{align}
\|(F_{\Phi},G_{\Phi},H_{\Phi})\|_{\dot{B}^{s}_{2,1}}\leq& C \|(a_{\Phi},u_{\Phi},\tau_{\Phi})\|_{\dot{B}^{n/2}_{2,1}}\|(\na a_{\Phi},\na u_{\Phi},\na \tau_{\Phi})\|_{\dot{B}^{s}_{2,1}}\nn\\
   \leq& C \|(a_{\Phi},u_{\Phi},\tau_{\Phi})\|_{\dot{B}^{n/2}_{2,1}}\|( a_{\Phi}, u_{\Phi}, \tau_{\Phi})\|_{\dot{B}^{s+1}_{2,1}},\nn
\end{align}
which is \eqref{fproduct8}. \qed

\section{Proof of Theorem \ref{thr3}}\label{sec3}

Inserting \eqref{fproduct8} into \eqref{apriori}, then there exits constant $\e_0$ and $C$ such that if $\|a\|_{L^\i}\leq \e_0$, we can obtain the following a priori estimate for system \eqref{cobr3first} with $F,G$ and $H$ being given in \eqref{fgh}.
\be\label{apriori1}
\bali
&\|(a_{\Phi},u_{\Phi},\tau_{\Phi})\|_{\dot{B}^{n/2-1,n/2}}+\la \|(a_{\Phi},u_{\Phi},\tau_{\Phi})\|_{L^1_{t,\dot{\th}}(\dot{B}^{n/2,n/2+1})}+\|(a_{\Phi},u_{\Phi},\tau_{\Phi})\|_{L^1_{t}(\dot{B}^{n/2}_{2,1})}\\
\leq& C\lt(\|e^{2\la_0\La}(a_0,u_0,\tau_0)\|_{\dot{B}^{n/2-1,n/2}}+\int^t_0\|(a_{\Phi},u_{\Phi},\tau_{\Phi})(s)\|_{(\dot{B}^{n/2}_{2,1})}\|( a_{\Phi}, u_{\Phi}, \tau_{\Phi})(s)\|_{({B}^{n/2,n/2+1}_{2,1})}ds\rt).
\eali
\ee

First, we approximate \eqref{cobr3first} by a sequence of ordinary differential equations by the classical Friedrichs method, see \cite{CheminGP:2011ANNALS} or \cite[Chapter 10]{BahouriCD:2011SPRINGER} for instance.

Let $L^2_k(\bR^n)$ denote the set of $L^2(\bR^n)$ functions spectrally supported in the annulus
\bes
\mathcal{C}_k:=\{\xi\in\bR^n:k^{-1}\leq |\xi|\leq k\}.
\ees
Define $\dot{\mathbb{E}}_k: L^2\Rightarrow L^2_k$ be the Friedrichs projector by
\bes
\mathcal{F}\dot{\mathbb{E}}_k U(\xi):={\bf1}_{\mathcal{C}_k}\mathcal{F}U(\xi),\q \text{for all}\q \xi\in\bR^n.
\ees
Consider the following ODE approximate system
\be\label{global1}
\bali
{\dfrac{d}{dt}}&\begin{pmatrix} a\\ u \\ \tau\end{pmatrix}=\dot{\mathbb{E}}_k
\begin{pmatrix} -\operatorname{div}-\na\cdot(au)\\ -\na a+\text{div}\tau -u\cdot \na u+[1-(1+a)^{\g-2}]\na a-\f{a}{1+a}\operatorname{div}\tau \\ -\tau+2\bD(u)-u\cdot\na \tau-Q(\tau,\na u)\end{pmatrix},
\\
&\begin{pmatrix} a\\ u \\ \tau\end{pmatrix}_{|t=0}=\dot{\mathbb{E}}_k\begin{pmatrix} a_0\\ u_0 \\ \tau_0\end{pmatrix}.
\eali
\ee
Solutions of \eqref{global1} is represented by $(a^k,u^k,\tau^k)$. Define the solution space by
\bes
\t{L}^2_k(\bR^n):=\{(a^k,u^k,\tau^k)\big | \text{inf}_{x\in\bR^n} | a>-1\}.
\ees
Note that if $\|a_0\|_{\dot{B}^{n/2}_{2,1}}$ is small, then $1+\dot{\mathbb{E}}_k a_0$ is positive for large $k$. then the initial data of \eqref{global1} are in $\t{L}^2_k(\bR^n)$. Thanks to the low-frequences cut-off of the operator $\dot{\mathbb{E}}_k$, all the Sobolev norms are equivalent. For fixed $k$, solving the ODE system \eqref{global1}, there exists a time maximal existing time $T_k$ such that
\bes
(a^k,u^k,\tau^k)\in C^1([0,T_k);\t{L}^2_k(\bR^n)).
\ees
For the obtained solution $(a^k,u^k,\tau^k)$, we define $\th_k(t)$ being the solution of the following ODE problem.
\be\label{radius}
\dot{\th_k}(t)=\|(a^k_{\Phi_k},u^k_{\Phi_k},\tau^k_{\Phi_k})\|_{\dot{B}^{n/2}_{2,1}}, \q\text{with}\q \th_k(0)=0,\q \Phi_k(t,\xi)=(2\la_0-\la \th_k(t))|\xi|.
\ee
Since the Fourier transform of $(a^k,u^k,\tau^k)$ is compactly supported, the righthand of the above ODE is Lipschitz with respect to $\th_k$. Then \eqref{radius} has a unique solution on $[0,T_k)$.

Next noting that the initial data $(a_0,u_0,\tau_0)$ of system \ref{cobr3first} is analytic  with $e^{2\la_0\La}(a_0,u_0,\tau_0)$ lying $\dot{B}^{\f{n}{2}-1,\f{n}{2}}$. It is also obviously that
\bes
\|e^{2\la_0\La}(a^k_0,u^k_0,\tau^k_0)\|_{\dot{B}^{\f{n}{2}-1,\f{n}{2}}}\leq \|e^{2\la_0\La}(a_0,u_0,\tau_0)\|_{\dot{B}^{\f{n}{2}-1,\f{n}{2}}}\leq \e_0\la_0.
\ees
Now we define $T^\ast_k$ to be
\bes
T^\ast_k:=\sup\{t\in [0,T_k): \th_k(t)\leq \f{\la_0}{\la},\q \text{and}\q \|(a^k,u^k,\tau^k)\|_{\dot{B}^{n/2-1,n/2}_{2,1}}\leq M\e_0\},
\ees
where $\la$ and $M$ are two constants to be determined later.

Next, we will show that by choosing suitable large $\la$ and $M$, we can obtain that $T^\ast_k=T_k$. We will use continuity argument and \eqref{apriori1} to show this.

Performing energy estimates almost the same as \eqref{apriori1} to system \eqref{global1}, we can obtain that for any $t\in [0,T^\ast_k)$, we have
\be\label{global2}
\bali
&\|(a^k_{\Phi},u^k_{\Phi},\tau^k_{\Phi})(t)\|_{\dot{B}^{n/2-1,n/2}}+\la \|(a^k_{\Phi},u^k_{\Phi},\tau^k_{\Phi})\|_{L^1_{t,\dot{\th}}(\dot{B}^{n/2,n/2+1})}+\|(a^k_{\Phi},u^k_{\Phi},\tau^k_{\Phi})\|_{L^1_{t}(\dot{B}^{n/2})}\\
\leq& C\lt(\|e^{2\la_0\La}(a^k_0,u^k_0,\tau^k_0)\|_{\dot{B}^{n/2-1,n/2}}+\int^t_0\|(a^k_{\Phi},u^k_{\Phi},\tau^k_{\Phi})(s)\|_{(\dot{B}^{n/2}_{2,1})}\|( a^k_{\Phi}, u^k_{\Phi}, \tau^k_{\Phi})(s)\|_{({B}^{n/2,n/2+1}_{2,1})}ds\rt)\\
\leq& C_0 \e_0\la_0+C_0\int^t_0\dot{\th}(s)\|( a^k_{\Phi}, u^k_{\Phi}, \tau^k_{\Phi})(s)\|_{({B}^{n/2,n/2+1}_{2,1})}ds.
\eali
\ee

Now, we first choose $\la=2C_0$, then from \eqref{global2}, we have
\be\label{global3}
\bali
&\|(a^k_{\Phi},u^k_{\Phi},\tau^k_{\Phi})\|_{\dot{B}^{n/2-1,n/2}}+C_0 \|(a^k_{\Phi},u^k_{\Phi},\tau^k_{\Phi})\|_{L^1_{t,\dot{\th}}(\dot{B}^{n/2,n/2+1})}+\|(a^k_{\Phi},u^k_{\Phi},\tau^k_{\Phi})\|_{L^1_{t}(\dot{B}^{n/2})}\\
\leq& C_0 \e_0\la_0.
\eali
\ee
Now, we choose $M:=2C_0\la_0$ and $\e_0$ sufficiently small such that
\bes
C_0\e_0\leq \f{1}{2C_0}=\f{1}{2\la}.
\ees
Then from \eqref{global3}, we can obtain that for any $t\in [0,T^\ast_k)$,
\bes
\th_k(t)\leq \f{\la_0}{2\la},\q \text{and}\q \|(a^k,u^k,\tau^k)\|_{\dot{B}^{n/2-1,n/2}_{2,1}}\leq \f{M}{2}\e_0.
\ees
By continuity, we can see that $T^\ast_k=T_k$ and the estimate \eqref{global3} is valid for any $t\in [0,T_k)$.  Also by the same continuity argument, we can see that system \eqref{global1} have a global time solution which satisfies for any $t\in [0,+\i)$, the estimate \eqref{global3} stand and for any $t\in [0,+\i)$,
\bes
\th_k(t)\leq \f{\la_0}{\la}.
\ees
Thanks to the uniform bound in \eqref{global3} and the uniformly low bound for the analytic radius $2\la_0-\la\th(t)\geq \la_0$, one can deduce, by a compactness argument, that there exists a unique solution $(a,u,\tau)$ to system \eqref{cobr3first} with the same bound in \eqref{global3}. The details are omitted. See for example \cite[Chapter 10]{BahouriCD:2011SPRINGER}. \qed

\vskip 1cm

\begin{appendix}

\section{Proof of Lemma \ref{lemproductf}}

We introduce the Bony decomposition \cite{Bony:1981ASENS} to perform nonlinear estimates in Besov spaces. The paraproduct between $f$ and $g$ is defined by
\bes
T_f g:=\sum_{k\in\bZ} \dot{S}_{k-1}f \dD_k g,
\ees
and the remainder is given by
\bes
R(f,g):=\sum_{|k-j|\leq 1} \dD_k f \dD_j g.
\ees
Then for $f,\, g\in\mathcal{Z}'(\bR^3) $, we have $fg=T_f g+T_g f+R(f,g)$.

We note that $\Phi(t,\xi)$ satisfies $\Phi(t,\xi)\leq \Phi(t,\eta)+\Phi(t,\xi-\eta)$ for $\xi,\eta\in\bR^3.$ Using \cite[Lemma 6.1]{Zi:2021AIHP}, we have the following paraproduct estimate.
\begin{itemize}
\item If $s\in\bR$, then there exists a positive constant $C$, depending on $s$, such that for $(f_{\Phi},g_{\Phi})\in L^\i\times\dot{B}^s_{2,1}$, we have
\be\label{paraproduct1}
      \|(T_{f}g)_{\Phi}\|_{\dot{B}^s_{2,1}}\leq C\|f_{\Phi}\|_{L^\i}\|g_{\Phi}\|_{\dot{B}^s_{2,1}};
\ee
\item If $s_1\in\bR$ and $s_2<0$, then there exists a positive constant $C$, depending on $s_1$ and $s_2$, such that for $(f_{\Phi},g_{\Phi})\in \dot{B}^{s_2}_{\i,\i}\times\dot{B}^{s_1}_{2,1}$, we have
\be\label{paraproduct2}
      \|(T_{f}g)_{\Phi}\|_{\dot{B}^{s_1+s_2}_{2,1}}\leq C\|f_{\Phi}\|_{\dot{B}^{s_2}_{\i,\i}}\|g_{\Phi}\|_{\dot{B}^{s_1}_{2,1}};
\ee
\item If $s_1+s_2>0$, then there exists a positive constant $C$, depending on $s_1$ and $s_2$, such that for $(f_{\Phi},g_{\Phi})\in \dot{B}^{s_1}_{2,\i}\times\dot{B}^{s_2}_{2,1}$, we have
\be\label{paraproduct3}
      \|(R(f,g))_{\Phi}\|_{\dot{B}^{s_1+s_2-n/2}_{2,1}}\leq C\|f_{\Phi}\|_{\dot{B}^{s_1}_{2,\i}}\|g_{\Phi}\|_{\dot{B}^{s_2}_{2,1}}.
\ee
\end{itemize}

Now we come to prove \eqref{fproduct1}. By using the interpolation $\dot{B}^{n/p}_{p,1}\hookrightarrow L^\i$ for $p\in [1,\i]$ and \eqref{paraproduct1}, we see that
\begin{align}
 \|(T_{f}g)_{\Phi}\|_{\dot{B}^{s}_{2,1}}\leq C  \|f_{\Phi}\|_{L^\i}\|g_{\Phi}\|_{\dot{B}^{s}_{2,1}} \leq  C\|f_{\Phi}\|_{\dot{B}^{n/2}_{2,1}}\|g_{\Phi}\|_{\dot{B}^{s}_{2,1}}.\label{fproduct2}
\end{align}
When $s=n/2$, the same as \eqref{fproduct2}, we have
\begin{align}
 \|(T_{g}f)_{\Phi}\|_{\dot{B}^{n/2}_{2,1}}\leq C  \|g_{\Phi}\|_{L^\i}\|f_{\Phi}\|_{\dot{B}^{n/2}_{2,1}} \leq  C\|f_{\Phi}\|_{\dot{B}^{n/2}_{2,1}}\|g_{\Phi}\|_{\dot{B}^{s}_{2,1}}.\nn
\end{align}
When $s<n/2$, by choosing $s_1=n/2$ and $s_2=s-n/2$ in \eqref{paraproduct2} and the interpolation $\dot{B}^{s}_{2,1}\hookrightarrow \dot{B}^{s-n/2}_{\i,\i}$, we see that
\begin{align}
 \|(T_{g}f)_{\Phi}\|_{\dot{B}^{s}_{2,1}}\leq& C  \|g_{\Phi}\|_{\dot{B}^{s-n/2}_{\i,\i}}\|f_{\Phi}\|_{\dot{B}^{n/2}_{2,1}}\leq C\|f_{\Phi}\|_{\dot{B}^{n/2}_{2,1}}\|g_{\Phi}\|_{\dot{B}^{s}_{2,1}}.\nn
\end{align}
By choosing $s_1=n/2$ and $s_2=s$ in \eqref{paraproduct3} and the interpolation $\dot{B}^{3/p}_{p,1}\hookrightarrow \dot{B}^{3/p}_{p,\i}$, we see that
\begin{align}
 \|(R(f,g))_{\Phi}\|_{\dot{B}^{s}_{2,1}}\leq& C \|f_{\Phi}\|_{\dot{B}^{n/2}_{2,\i}}\|g_{\Phi}\|_{\dot{B}^{s}_{2,1}}\leq C \|f_{\Phi}\|_{\dot{B}^{n/2}_{2,1}}\|g_{\Phi}\|_{\dot{B}^{s}_{2,1}}. \label{fproduct5}
\end{align}

Combining estimates in \eqref{fproduct2} to \eqref{fproduct5}, we obtain \eqref{fproduct1}. \qed

\end{appendix}

\section*{Data availability statement}

\q\ Data sharing is not applicable to this article as no datasets were generated or analysed during the current study.

\section*{Conflict of interest statement}

\q\ The authors declare that they have no conflict of interest.

\section*{Acknowledgments}
\q  X. Pan thanks to Professor Ruizhao Zi for helpful discussion and this project is supported by National Natural Science Foundation of China (No. 11801268, 12031006).
\bibliographystyle{plain}

%%%%%%%%%%%%%%%%%%%%%%%%%%%%%%%%%%%%

\end{document}